\documentclass[12pt]{amsart}
\usepackage[margin=2.0cm]{geometry}

\newenvironment{ack}{\textbf{Acknowledgments.}}{\par}
\newenvironment{funding}{}{}

\usepackage{amssymb}
\usepackage{upgreek}
\usepackage{longtable}
\usepackage{tabularx}
\usepackage{hhline}
\usepackage{boldline,multirow}
\usepackage[matrix,arrow,curve]{xy}
\usepackage{extarrows }
\usepackage{calligra,mathrsfs}

\DeclareMathOperator{\Hom}{\mathscr{H}\text{\kern -3pt {\calligra{om}}}\,}
\newcommand{\Supp}{\operatorname{Supp}}

\newcommand{\Exc}{\operatorname{Exc}}

\newcommand{\p}{\mathrm{p}_{\mathrm{a}}}
\newcommand{\Sing}{\operatorname{Sing}}
\newcommand{\red}{\operatorname{red}}

\newcommand{\Pic}{\operatorname{Pic}}

\newcommand{\CC}{\mathbb{C}}
\newcommand{\ZZ}{\mathbb{Z}}
\newcommand{\PP}{\mathbb{P}}
\newcommand{\QQ}{\mathbb{Q}}

\newcommand{\LLL}{{\mathscr{L}}}
\newcommand{\OOO}{{\mathscr{O}}} 
\newcommand{\mumu}{{\boldsymbol{\mu}}}

\def\typec#1{$\mathrm{{\left(#1\right)}}$}
\def\type#1{$\mathrm{{#1}}$}
\def\mtypec#1{\mathrm{{\left(#1\right)}}}
\def\mtype#1{\mathrm{{#1}}}

 \renewcommand{\MR}[1]{}

\newcounter{NN}

\newtheorem{theorem}[subsection]{Theorem}
\newtheorem{proposition}[subsection]{Proposition}

\newtheorem{corollary}[subsection]{Corollary}
\newtheorem{conjecture}[subsection]{Conjecture}

\newtheorem{principle-conjecture}[subsection]{Principle-Conjecture}
\newtheorem*{plain*}{}

\theoremstyle{definition}
\newtheorem*{definition*}{Definition}
\newtheorem{definition}[subsection]{Definition}

\newtheorem{example}[subsection]{Example}
\newtheorem*{example*}{Example}
\newtheorem*{examples*}{Examples}
\newtheorem{example-construction}[subsection]{Example-Construction}

\newtheorem*{construction*}{Construction}

\newtheorem*{notation*}{Notation}

\numberwithin{equation}{section}

\begin{document}

 \makeatletter
 \@addtoreset{equation}{subsection}
 \makeatother

\renewcommand{\thesubsection}{\arabic{section}.\arabic{subsection}}
\renewcommand{\theequation}{\arabic{section}.\arabic{subsection}.\arabic{equation}}


\title{Effective results in the three-dimensional minimal model program}

\author{Yuri Prokhorov}
\address{Steklov Mathematical Institute,
8 Gubkina street, Moscow 119991, Russia }
\email{prokhoro@mi-ras.ru}


\begin{abstract}
We give a brief review on recent developments in 
the three-dimensional minimal model program.
\end{abstract}

\maketitle


In this note we give a brief review on recent developments in 
the three-dimensional minimal model program (MMP for short).
Certainly, this is not a complete survey of all advances in this area. 
For example, we do not discuss the minimal models of
varieties of non-negative Kodaira dimension, as well as, applications to birational geometry and moduli spaces.

The aim of the MMP is to find a good representative 
in a fixed birational equivalence class of algebraic varieties.
Starting with an arbitrary smooth projective variety one can perform a finite number of certain elementary transformations, called
divisorial contractions and flips, and at the end obtain a variety which is 
simpler in some sense. Most parts of the MMP are completed in arbitrary dimension. One of the basic remaining problems is the following: 
\par\smallskip
\emph{Describe all the intermediate steps 
and the outcome of the MMP.}
\par\smallskip
The MMP makes sense only in dimensions $\ge 2$ 
and for surfaces it is classical and well-known. 
So the first non-trivial case is the three-dimensional one.
It turns out that to proceed with the MMP in dimension $\ge 3$ one has to work 
with varieties admitting certain types of very mild, so-called 
terminal, singularities.
On the other hand, dimension $3$ is the last dimension where one can expect effective results: in higher dimensions classification results become very complicated and unreasonably long.

We will work over the field $\CC$ of complex numbers throughout.
A variety is either an algebraic variety or a reduced complex space.

\section{Singularities}
Recall that a Weil divisor $D$ on a normal variety is said to be \emph{$\QQ$-Cartier} if its multiple $nD$, for some $n$, is a Cartier divisor.
For any morphism $f: Y\to X$, the pull-back $f^*D$ of a $\QQ$-Cartier divisor $D$ is well defined as a divisor with rational coefficients ($\QQ$-divisor).
A variety $X$ has \emph{$\QQ$-factorial singularities} if any 
Weil divisor on $X$ is $\QQ$-Cartier.

\begin{definition}
\label{def:sing}
A normal algebraic variety (or an analytic space) $X$ is said to have \emph{terminal} (resp. \emph{canonical}, \emph{log terminal}, \emph{log canonical}) singularities if the canonical Weil divisor $K_X$ is $\QQ$-Cartier and
for any birational morphism
$f: Y\to X$ one can write
\begin{equation}
\label{eq:discr}
K_Y =f^*K_X+\sum a_i E_i, 
\end{equation} 
where $E_i$ are all the exceptional divisors and $a_i>0$ (resp. $a_i\ge 0$, $a_i>-1$, $a_i\ge -1$) for all $i$. The smallest positive $m$ such that $mK_X$ is Cartier is called the \emph{Gorenstein index} of $X$.
Canonical singularities of index $1$ are rational Gorenstein.
\end{definition}
The class of terminal $\QQ$-factorial singularities is the smallest class that 
is closed under the MMP. Canonical singularities important because they appear on the canonical models 
of varieties of general type.
A crucial observation is that terminal 
singularities lie in codimension $\ge 
3$. In particular, terminal surface singularities are smooth
and terminal threefold singularities are isolated.
Canonical singularities of surfaces are called \emph{Du Val} or rational double points. Any two-dimensional log terminal singularity is a quotient of a smooth surface germ \cite{Kawamata:crep}.
Terminal threefolds singularities were classified by M.~Reid \cite{Reid:MM} and S.~Mori \cite{Mori:term-sing}.

\begin{example*}
Let $X\subset \CC^4$ be a hypersurface given by the equation
\[
\phi(x_1,x_2,x_3)+x_4\psi(x_1,\dots,x_4) =0,
\]
where $\phi=0$ is an equation of a Du Val (ADE) singularity. Then the singularity of $X$ at $0$ is canonical Gorenstein. 
It is terminal if and only if it is isolated.
Singularities of this type are called \type{cDV}. 
\end{example*}

According to \cite{Reid:MM} any three-dimensional terminal singularity of index $m>1$ is a quotient of an isolated \type{cDV}-singularity by
the cyclic group $\mumu_m$ of order $m$. More precisely, we have the following

\begin{theorem}[{\cite{Reid:MM}}]
Let $(X\ni P)$ be an analytic germ
of a three-dimensional terminal singularity of index $m\ge 1$. Then there exists an isolated \type{cDV}-singularity $\big(X^\sharp\ni P^\sharp\big )$
and a cyclic $\mumu_m$-cover
\[
\pi: \big(X^\sharp\ni P^\sharp \big)\longrightarrow (X\ni P )
\]
which is \'etale outside $P$. 
\end{theorem}

The morphism $\pi$ in the above theorem is called the \emph{index-one cover}.
A detailed classification of all possibilities for the equations of $X^\sharp \subset \CC^4$ and the actions of $\mumu_m$ was obtained in \cite{Mori:term-sing}
(see also \cite{Reid:YPG}). 

\begin{example*}
Let the cyclic group $\mumu_m$ act on $\CC^n$ diagonally via
\[
(x_1,\dots, x_n) \longmapsto \big(\upzeta^{a_1} x_1,\, \dots, \upzeta^{a_n}x_n\big),\quad \upzeta=\upzeta_m=\exp (2\uppi \operatorname{i}/m).
\]
Then we say that $(a_1,\dots,a_n)$ is the collection of \emph{weights} of the action.
Assume that the action is free in codimension $1$.
Then the quotient singularity $\CC^n/\mumu_m\ni 0$ is said to be of type $\frac1m(a_1,\dots,a_n)$.
According to the criterion (see \cite[Theorem~4.11]{Reid:YPG}) this singularity 
is terminal if and only if 
\[
\sum_{i=1}^{n} \overline{k a_i}> m\text{\quad for $k=1,\dots, m-1$,} 
\]
where $\overline{\phantom{a_i}}$ is the smallest residue $\mod m$. In the threefold case this criterion has a very simple form: 
a quotient singularity $\CC^m/\mumu_m$ is terminal if and only if 
it is of type $\frac1{m}(1,-1,a)$, where $\gcd(m,a)=1$. 
This is a 
\emph{cyclic quotient terminal singularity}.
\end{example*}

\begin{example*}[\cite{Mori:term-sing, Reid:YPG}]
Let the cyclic group $\mumu_m$ act on $\CC^4$ diagonally 
with weights $(1,-1,a,0)$, where $\gcd(m,a)=1$. 
Then for a polynomial $\phi(u,v)$ the singularity at~$0$ of the quotient 
\[
\left\{x_1x_2+\phi(x_3^m,x_4)=0\right\}/\mumu_m
\]
is
terminal whenever it is isolated. 
The index of this singularity equals $m$.
\end{example*}

As a consequence of the classification, we obtain that the local fundamental group of the (analytic) germ of a three-dimensional terminal singularity of index $m$ is cyclic of order~$m$:
\begin{equation}
\label{eq:pi1}
\uppi_1 (X\setminus \{P \}) \simeq \ZZ/m\ZZ.
\end{equation} 
In particular, for any Weil $\QQ$-Cartier divisor $D$ on $X$ its $m$th multiple $mD$ is Cartier \cite[Lemma~5.1]{Kawamata:crep}. 

The class of canonical threefold singularities is much larger than the class of terminal ones. However there are certain boundedness results. For example, it is known that  
the index of a strictly canonical isolated singularity is at most $6$ \cite{Kawakita-index}.

The modern higher dimensional MMP often works with pairs and one need to extend Definition~\ref{def:sing} to a wider class of objects.
\begin{definition*}
Let $X$ be a normal variety and let $B$ be an effective $\QQ$-divisor on $X$. The pair $(X,B)$ is said to be \emph{plt} (resp. \emph{lc}) if 
$K_X+B$ is $\QQ$-Cartier and
for any birational morphism
$f: Y\to X$ one can write
\[
K_Y +B_Y=f^*(K_X+B)+\sum a_i E_i, 
\]
where $B_Y$ is the proper transform of $B$, \ $E_i$ are all the exceptional divisors and $a_i>-1$ (resp. $a_i\ge -1$) for all $i$. The pair $(X,B)$ is said to be \emph{klt} if it is plt and $\lfloor B\rfloor=0$.
\end{definition*}

\section{Minimal Model Program}
Basic elementary operations in the MMP are Mori contractions.

A \emph{contraction} is a proper surjective morphism $f:X\to Z$ of normal 
varieties 
with connected fibers. The \emph{exceptional locus} of a contraction $f$ is the subset $\Exc(f)\subset X$ of points 
at which $f$ is not an isomorphism. 
A \emph{Mori contraction} is a contraction $f:X\to Z$
such that 
the variety $X$ has at worst terminal 
$\QQ$-factorial singularities, the anticanonical class ${-}K_X$ is $f$-ample, and the relative Picard number $\uprho(X/Z)$ equals $1$.
A Mori contraction is said to be \emph{divisorial} (resp. \emph{flipping}) if it is birational and the locus
$\Exc(f)$ has codimension $1$ (resp. $\ge 2$).
For a divisorial contraction the exceptional locus $\Exc(f)$ is a prime divisor. 
A Mori contraction, whose target is a lower dimensional variety, is called \emph{Mori fiber space}.
Then the general fiber is a Fano variety with at worst terminal singularities. In the particular cases where the relative dimension 
of $X/Z$ equals $1$ (resp. $2$) the Mori fiber space $f:X\to Z$ 
is called \emph{$\QQ$-conic bundle} (resp. \emph{$\QQ$-del Pezzo fibration}).
If $Z$ is a point, then $X$ is a Fano variety with at worst terminal $\QQ$-factorial singularities and $\Pic(X)\simeq \ZZ$. For short, we call such 
varieties \emph{$\QQ$-Fano}.

The MMP procedure is a sequence of elementary transformations which are constructed inductively \cite{KMM,KM:book}. Let $X$ be a projective algebraic variety 
with terminal $\QQ$-factorial singularities.
If the canonical divisor $K_X$ is not nef, then there exists a Mori contraction $f:X\to Z$. 
If $f$ is divisorial, then 
$Z$ is again a variety with terminal $\QQ$-factorial singularities
and, in this situation, we can proceed with the MMP replacing $X$ with $Z$.
In contrast, a flipping contraction takes us out the category of terminal $\QQ$-factorial varieties. To proceed, one has to perform a surgery operation as follows
\[
\xymatrix@R=0.8em{
X\ar[dr]_f\ar@{-->}[rr]&&X^+\ar[dl]^{f^+}
\\
&Z&
} 
\]
where $f^+$ is a contraction whose exceptional locus has codimension $\ge 2$ and the divisor
$K_{X^+}$ is $\QQ$-Cartier and $f^+$-ample. Then the variety $X^+$ again has 
terminal $\QQ$-factorial singularities and we can proceed by replacing $X$ with $X^+$.

The process described above should terminate and at the end we obtain a variety $\bar X$ such that either $\bar X$ has a Mori fiber space structure $\bar X\to \bar Z$ or 
$K_{\bar X}$ is nef. 
One of the remaining open problems is the termination of the program, to be more precise, 
termination of a sequence of flips. The problem solved affirmatively
in dimension $\le 4$ \cite{Shokurov:non-van,KMM}, for varieties of general type, for uniruled varieties \cite{BCHM}, and in some other special cases.
We refer to \cite{Birkar:ICM} for more comprehensive survey of the higher-dimensional MMP.

The MMP has a huge number of applications in algebraic geometry.
The most impressive consequence of the MMP is the finite generation of the canonical ring 
\[
\mathrm{R}(X,K_X) := \bigoplus_{n\ge 0} H^0(X,\OOO_X(mK_X))
\]
of a smooth projective variety $X$ \cite{BCHM,Hacon-McKernan:R2}.
Another application of the MMP is so-called Sarkisov program which 
allows to decompose a birational maps between Mori fiber spaces into 
composition of elementary transformations, called Sarkisov links \cite{Corti95:Sark, HaconMcKernan:Sark,Shokurov-Choi}.
Also the MMP can be applied to varieties with finite group actions and 
to varieties over non-closed fields (see~\cite{P:G-MMP}).

As was explained above, the Mori contractions are fundamental building blocks in the MMP. To apply the MMP effectively, one needs to understand the structure of its steps in details.
For a Mori contraction $f:X\to Z$ of a three-dimensional variety $X$ there are only the following possibilities:
\begin{itemize}
\item 
$f$ is divisorial and the image of the (prime) exceptional divisor $E:=\Exc(f)$
is either a point or an irreducible curve,

\item
$f$ is flipping and the exceptional locus $\Exc(f)$ is a union of a finite number of irreducible curves,

\item
$Z$ is a surface and $f$ is a $\QQ$-conic bundle,

\item
$Z$ is a curve and $f$ is a $\QQ$-del Pezzo fibration,

\item
$Z$ is a point and $X$ is a $\QQ$-Fano threefold.
\end{itemize} 
Mori contractions of smooth threefolds to varieties of positive dimension where classified in the pioneer work of S.~Mori \cite{Mori:3-folds}.
S.~Cutkosky \cite{Cutkosky:contr} extended this 
classification to the case of Gorenstein terminal varieties. Smooth Fano threefolds of Picard number one where classified by Iskovskikh \cite{Isk:Fano1e,Isk:Fano2e} (see also \cite{IP99}).
 Fano threefolds with Gorenstein terminal 
singularities are degenerations of smooth ones \cite{Namikawa:Fano}.
Below we are going to discuss Mori contractions of threefolds. We are interested only in the biregular structure of a contraction $f:X\to Z$ near a fixed fiber $f^{-1}(o)$, $o\in Z$. Typically we do not consider the simple case where $X$ is Gorenstein.

\section{General elephant}
\label{sect:ge}
A natural way to study higher-dimensional varieties is the inductive one. 
Typically to apply this method we need to find a certain subvariety of dimension one less (divisor)
which is sufficiently good is the sense of singularities.
\begin{conjecture}
\label{conj:ge}
Let $f: X\to (Z\ni o)$ be a threefold Mori contraction, 
where $(Z\ni o)$ is a small neighborhood. 
Then the general member $D\in |{-}K_X|$ is a normal surface with Du Val 
singularities.
\end{conjecture}
The conjecture was proposed by M.~Reid who called a good member of $|{-}K_X|$ ``elephant''. We follow this language and call \ref{conj:ge} the
General Elephant Conjecture.
The importance of the existence of good member in $|{-}K_X|$ is motivated by many reasons:
\begin{itemize}
 \item 
The general elephant passes through all the non-Gorenstein points of $X$ and so it encodes the information about their types and configuration (cf. Proposition~\ref{prop:local-ge} below). 
\item 
For flipping contractions Conjecture~\ref{conj:ge} is a sufficient condition 
for the existence of threefold flips \cite{Kawamata:crep}.

\item 
For a divisorial contraction $f: X\to Z$ whose fibers have dimension $\le 1$ the image $D_Z:=f(D)$ of a Du Val elephant $D\in |{-}K_X|$ must be again Du Val and the image $\Gamma:=f(E)$ of the exceptional divisor is a curve on $D_Z$. Then one can reconstruct $f$ starting from the triple $(Z\supset D_Z\supset \Gamma)$ by using certain birational procedure.
Such an approach was successfully worked out in many cases by N.~Tziolas
\cite{Tziolas:E6-E7, Tziolas:CompMath03, Tziolas:ss, Tziolas:D-mm+err}.

 \item 
If $f: X\to (Z\ni o)$ is a $\QQ$-del Pezzo fibration such that general $D\in |{-}K_X|$ is Du Val, then composing the projection $D\to Z$ with 
minimal resolution $\tilde D\to D$ we obtain a relatively minimal 
elliptic fibration whose singular fibers are classified by Kodaira \cite{Kodaira:surfaces2-3}. Then one can get a bound of multiplicities of fibers and describe the configuration of non-Gorenstein singularities.
 \item 
For a $\QQ$-Fano threefold $X$, a Du Val general elephant is a (singular) K3 surface. In the case where the linear system $|{-}K_X|$ is ``sufficiently big''
this implies the existence of a good Gorenstein model \cite{Alexeev:ge}.
\end{itemize}

Shokurov \cite{Shokurov:flips+err} generalized Conjecture~\ref{conj:ge} and introduced a new notion which is very efficient in the study of pluri-anticanonical linear systems. Omitting technicalities we reproduce a weak form of Shokurov's definition. 

\begin{definition*}
An \emph{$n$-complement} of 
the canonical class 
$K_X$ is a member $D\in |{-}nK_X|$ such that 
the pair $(X,\frac 1n D)$ is lc. 
An $n$-complement is said to be klt (resp. plt) if so 
the pair $(X,\frac 1n D)$~is. 
\end{definition*}
According to the inversion of adjunction \cite{Shokurov:flips+err} the existence of 
a Du Val general elephant $D\in |{-}K_X|$ is equivalent to the existence 
of a plt $1$-complement. Shokurov developed a powerful theory that works in arbitrary dimension and allows to
construct complements inductively (see \cite{Shokurov:flips+err}, \cite{P:lect-compl} and references therein).

Note that 
Reid's general elephant fails for Fano threefolds. For example, in \cite{GRD, Fletcher:wci} one can find examples of $\QQ$-Fano threefolds 
with empty anticanonical linear system. 
Because of this, the statement of~\ref{conj:ge} sometimes is called ``principle''.
Nonetheless there are only a few examples of such Fano threefolds.  In the cases $\dim(Z)>0$
Conjecture~\ref{conj:ge} is expected to be true. 
The following should be considered as the local version of~\ref{conj:ge}.

\begin{proposition}[Reid {\cite{Reid:YPG}}]
\label{prop:local-ge}
Let $(X\ni P)$ be the analytic germ of a threefold terminal singularity of index $m>1$. Then the general 
member $D\in |{-}K_X|$ is a Du Val singularity.
Furthermore, let $\pi: X'\to X$ be the index-one cover and let $D':=\pi^{-1}(D)$. 
Then the cover $D' \to D$ belongs to one of 
the following six types:
\begin{center}
{\renewcommand\arraystretch{1.1}
\begin{tabularx}{\textwidth}{lp{8em}@{\hspace{7em}}|lp{7em} }
\textup{$(X\ni P)$} & $D' \longrightarrow D$ &\textup{$(X\ni P)$} & $D' \longrightarrow D$
\\[-0.7em]
\multicolumn{4}{l}{\rule{0.9\textwidth}{.3mm}}
\\[-0.4em]
\type{cA/m} & $\mtype{A_{k-1}} \xrightarrow{\ m:1\ } \mtype{A_{km-1}}$ & \type{cAx/2} & $\mtype{A_{2k-1}} \xrightarrow{\ 2:1\ } \mtype{D_{k+2}}$
\\ \hline
 \type{cAx/4} & $\mtype{A_{2k-2}} \xrightarrow{\ 4:1\ } \mtype{D_{2k+1}}$& 
\type{cD/2} & $\mtype{D_{k+1}} \xrightarrow{\ 2:1\ } \mtype{D_{2k}}$
\\ \hline
 \type{cD/3} & $\mtype{D_{4}} \xrightarrow{\ 3:1\ } \mtype{E_{6}}$& 
\type{cE/2} & $\mtype{E_{6}} \xrightarrow{\ 2:1\ } \mtype{E_{7}}$
\end{tabularx}}
\end{center}
\end{proposition}

\section{Divisorial contractions to a point}
In this section we treat divisorial Mori contractions of a divisor to a point.
This kind of contractions is studied very well
due to works of
Y.~Kawamata \cite{Kawamata:Div-contr}, A.~Corti \cite{Corti2000}, M.~Kawakita \cite{Kawakita:smooth,Kawakita:A1,Kawakita:ge, Kawakita:hi,Kawakita:supp}, T.~Hayakawa \cite{Hayakawa:blowup1,Hayakawa:blowup2,Hayakawa:discr1}, and others.
In this case General Elephant Conjecture~\ref{conj:ge} has been proved:

\begin{theorem}[Kawakita \cite{Kawakita:ge, Kawakita:hi}]
Let $f: X\to (Z\ni o)$ be a divisorial Mori contraction that contracts a 
divisor to a point.
Then the general member $D\in |{-}K_X|$ is Du Val. 
\end{theorem}
One of the main tools in the proofs is the orbifold Riemann-Roch formula \cite{Reid:YPG}: if $X$ is a three-dimensional projective variety with terminal singularities and $D$ is a Weil $\QQ$-Cartier divisor on $X$, then for the sheaf $\LLL=\OOO_X(D)$ 
there is a formula of the form
\begin{equation}
\label{eq:RR}
\chi(\LLL)=\chi(\OOO_X)+\frac{1}{12}D\cdot (D-K_X)\cdot (2D-K_X) +\frac{1}{12} D\cdot \mathrm{c}_2 +\sum_{P} c_P(D), 
\end{equation} 
where the sum rungs through all the virtual quotient singularities of $X$, i.e.
through the actual singularities of $X$ are replaced with their small deformations \cite{Reid:YPG}, and $c_P(D)$ is a local contribution due to singularity at $P$,
depending only on the local analytic type of $D$ at $P$. There is an explicit formula for computation of $c_P(D)$.

Except for a few hard cases the classification of divisorial Mori contractions of a divisor to a point has been completed. 
Typical result here is to show that a contraction is a weighted blowup 
with some explicit collection of weights:

\begin{theorem}[Y. Kawamata {\cite{Kawamata:Div-contr}}]
Let $f: X\to (Z\ni o)$ be a divisorial Mori contraction that contracts a divisor to a 
point. Assume that $o\in Z$ is a cyclic quotient singularity of type 
$\frac 1r(a,-a,1)$. 
Then $f$ is the weighted blowup with weights $(a/r, 1 -a/r, 1/r)$. 
\end{theorem}

\begin{theorem}[M. Kawakita {\cite{Kawakita:smooth}}]
Let $f: X\to (Z\ni o)$ be a divisorial Mori contraction that contracts a divisor to a 
smooth point. 
Then $f$ is the weighted blowup with weights $(1, a, b)$, where $\gcd(a,\, 
b)=1$. 
\end{theorem}

These results are intensively used in the three-dimensional birational geometry, for example, in the proof of birational rigidity of index $1$ Fano threefold weighted hypersurfaces~\cite{Corti-Pukhlikov-Reid}.

\section{Del Pezzo fibrations}
Much less is known about local structure of $\QQ$-del Pezzo fibrations.
As was explained in Sect.~\ref{sect:ge}, the existence of a Du Val general elephant would be very helpful in the 
study these kind of contractions. 
However, in this case Conjecture~\ref{conj:ge} is established only in some special situation.

An important question that can be asked in the Del Pezzo fibration case is the presence of multiple fibers.

\begin{theorem}[{\cite{MP:DP-e}}]
\label{thm:DP-mult}
Let $f: X\to Z$ be a $\QQ$-del Pezzo fibration and let $f^*(o)=m_oF_o$
be a special fiber of multiplicity $m_o$.
Then $m_o\le 6$ and all the cases $1\le m_o\le 6$ occur.
Moreover, the possibilities for 
the local behavior of $F_o$ near singular points are described.
\end{theorem}
The main idea of the proof is to apply the orbifold Riemann-Roch formula \eqref{eq:RR} to the divisor  $F_o$ and its multiples.

\begin{example*}
Suppose that the cyclic group 
$\mumu_4$ acts on $\PP^1_{x}\times \PP^1_{y}\times \CC_t$
via 
\[
(x, y; t) \longmapsto \big(y,\ - x,\ \sqrt{-1}\ t\big). 
\]
Then the quotient
\[
X= \big(\PP^1\times \PP^1\times \CC\big)/\mumu_4\longrightarrow Z=\CC/\mumu_4
\]
is the germ of a $\QQ$-del Pezzo fibration with central fiber of multiplicity $4$.
\end{example*}

Another type $\QQ$-del Pezzo fibrations which is investigated relatively well are those whose 
central fiber $F:=f^{-1}(o)$ is reduced, normal, and has ``good'' singularities. 
Then $X$ can be viewed as a one-parameter \emph{smoothing} of $F$. The total space of this 
smoothing must be \emph{$\QQ$-Gorenstein} and $F$ can be viewed as a \emph{degeneration} of a general fiber (smooth del 
Pezzo surface) in a $\QQ$-Gorenstein family. 
The most studied class of singularities admitting $\QQ$-Gorenstein smoothings 
is the class of singularities of type~\type{T}.

\begin{definition*}[Koll\'ar, Shepherd-Barron {\cite{Kollar-ShB-1988}}]
A two-dimensional quotient singularity is said to be of \emph{type ~\type{T}} 
if it admits a smoothing in a one-parameter $\QQ$-Gorenstein family $X\to B$. 
\end{definition*}

In this case, by the inversion of adjunction \cite{Shokurov:flips+err}, the pair $(X, F)$ is plt and the total family $X$ is terminal.
Conversely, if $X\ni P$ is a $\QQ$-Gorenstein point and $F$ is an effective 
Cartier divisor at $P$ such that the pair $(X, F)$ is plt, then $F\ni P$ is a 
\type{T}-singularity and the point $X\ni P$ is terminal.
Singularities of type~\type{T} and their deformations were studied by Koll\'ar and Shepherd-Barron \cite{Kollar-ShB-1988}. In particular, they proved that 
any \type{T}-singularity is either a Du Val point or a 
cyclic quotient of type $\frac1{m}(q_1,q_2)$ with 
\[
\gcd(m,q_1)=\gcd(m,q_2)=1,\qquad (q_1+q_2)^2\equiv 0 \mod m.
\]
Minimal resolutions of these singularities are also described \cite[\S~3]{Kollar-ShB-1988}.

Thus to study $\QQ$-del Pezzo fibrations whose central fiber has only 
quotient singularities one has to consider $\QQ$-Gorenstein smoothings 
of del Pezzo surfaces with singularities of type~\type{T}.
The important auxiliary fact here is the unobstructedness of deformations:
\begin{proposition}[{\cite{Manetti:P2, HP10}}]
Let $F$ be a projective surface with log canonical singularities such that 
${-}K_F$
is big. Then there are no local-to-global obstructions to deformations of $F$. 
In particular, if $F$
has \type{T}-singularities, then $F$ admits a $\QQ$-Gorenstein smoothing.
\end{proposition}

\begin{theorem}[Hacking-Prokhorov {\cite{HP10}}]
\label{thm:DP-deg}
Let $F$ be a projective surface with quotient singularities such that ${-}K_F$ is
ample, $\uprho(F) = 1$, and $F$ admits a $\QQ$-Gorenstein smoothing. Then $F$ 
belongs to one of the following series:
\begin{itemize} 
\item 
$14$ infinite series of toric surfaces \textup(see below\textup);
\item 
partial smoothing of a toric surface as above;
\item
$18$ sporadic families of surfaces of index $\le 2$ \cite{Alexeev-Nikulin:MSJ}.
\end{itemize}
\end{theorem}
Toric surfaces appeared in the above theorem are determined by a Markov-type equation. More precisely, for $K_F^2\ge 5$ these surfaces are weighted projective spaces given by the following table:
\begin{center}
\begin{tabular}{@{\hspace{4em}}l@{\hspace{3em}}l@{\hspace{3em}}l}
$K_F^2$&\multicolumn{1}{c}{$F$}&Markov-type equation
\\[-0.7em]
\multicolumn{3}{l}{\rule{0.87\textwidth}{.3mm}}
\\[-0.4em]
9&$\PP\left(a^2, b^2, c^2\right)$ & $a^2 + b^2 + c^2 = 3abc $
\\
8&$\PP\left(a^2, b^2, 2c^2\right)$ &$a^2 + b^2 + 2c^2 = 4abc$
\\
6 &$\PP\left(a^2, 2b^2, 3c^2\right)$ & $a^2 + 2b^2 + 3c^2 = 6abc$
\\
5 &$\PP\left(a^2, b^2, 5c^2\right)$ & $a^2 + b^2 + 5c^2 = 5abc$
\end{tabular}
\end{center}
and for $K^2\le 4$ they are 
certain 
abelian quotients of the weighted projective spaces as above.
Note however that in general we cannot assert that for central fiber 
$F$ of a $\QQ$-del Pezzo fibration the condition $\uprho(F)=1$ holds.
Some partial results in the case $\uprho(F)>1$ where obtained in \cite{P:AIF15}. In particular, \cite{P:AIF15} establishes the existence of Du Val
general elephant for $\QQ$-del Pezzo fibrations with ``good'' fibers:
\begin{theorem}
Let $f: X\to (Z\ni o)$ be a $\QQ$-del Pezzo fibration over a curve germ $Z\ni o$. Assume that the fiber $f^{-1}(o)$ is reduced, normal and has only log terminal singularities. Then the
general elephant $D\in |{-}K_X|$ is Du Val. 
\end{theorem}

Theorem~\ref{thm:DP-deg} gives a complete answer to the question posed by M.~Manetti~\cite{Manetti:P2}:
\begin{corollary}[{\cite{HP10}}]
Let $X$ be a projective surface with quotient singularities which admits a
smoothing to $\PP^2$. Then $X$ is a $\QQ$-Gorenstein deformation of a weighted 
projective plane
$\PP\left(a^2, b^2, c^2 \right)$, where the triple $(a,b,c)$ is a solution of the Markov equation
\[
a^2+ b^2+ c^2=3a b c.
\]
\end{corollary}

Results similar to Theorem~\ref{thm:DP-deg} were obtained for $\QQ$-del Pezzo fibrations whose central fiber is log canonical \cite{P:lcdP:19}. However in this case the classification is not complete. 

\section{Extremal curve germs}
To study Mori contractions with fibers of dimension $\le 1$ it is convenient to work with analytic threefolds and 
to localize to situation near a curve contained in a fiber. 

\begin{definition}
Let $(X\supset C)$ be the \emph{analytic} germ 
of a threefold with terminal singularities along a reduced connected complete 
curve. Then
$(X\supset C)$ is called an \emph{extremal curve germ} if there exists a contraction 
\begin{equation*}
f: (X\supset C)\to (Z\ni o) 
\end{equation*}
such that\quad $C=f^{-1}(o)_{\red}$\quad and\quad ${-}K_X$ is $f$-ample. 
The curve $C$ is called the \emph{central fiber} of the germ and $Z\ni o$ is called the \emph{target variety} or the \emph{base} of $(X\supset C)$.
An extremal curve germ is said to be \emph{irreducible} if so its central fiber is.
\end{definition}

In the definition above we do not assume that $X$ is $\QQ$-factorial nor $\uprho(X/Z)=1$. This is because $\QQ$-factoriality typically
is not a local condition in the analytic category (see \cite[\S~1]{Kawamata:crep}).
There are three types of extremal curve germs.
\begin{itemize}
\item 
\emph{flipping} if is $f$ birational and does not contract divisors;
\item 
\emph{divisorial} if the exceptional locus is two-dimensional;
\item 
\emph{$\QQ$-conic bundle germ} if the target variety $Z$ is a surface.
\end{itemize}
If a divisorial curve germ is irreducible, then 
the exceptional locus of the corresponding contraction is a $\QQ$-Cartier divisor and the target variety $Z$ has terminal singularities \cite[\S 3]{MP:IA}.
In general this is not true. It may happen that the exceptional locus is a union of a divisor and some curves. 

As an example we consider that case where $X$ has 
singularities of indices $1$ and $2$.

\begin{theorem}[{\cite{MP:cb1}}]
\label{thm:cb-index2}
Let $(X\supset C)$ be a $\QQ$-conic bundle germ 
over a smooth base. 
Assume that $X$ is not Gorenstein and $2K_X$ is Cartier. Then 
$X$ can be embedded to $\PP(1,1,1,2)\times \CC^2$ and given there by two quadratic equations. In particular,
the point $P\in X$ of index $2$ is unique, the curve $C$ has at most $4$ components, all of them pass through 
$P$.
\end{theorem}

\begin{theorem}[{\cite{KM:92}}]
\label{thm:indexdf}
Let $(X\supset C)$ be a flipping extremal curve germ and let 
\[
\xymatrix@R=1em{
&(X\supset C)\ar@/_/[dr]^f
\ar@{-->}[rr] && (X^+\supset C^+)\ar@/^/[dl]_{f^+}
\\
&&(Z\ni o)&
}
\]
be the corresponding flip. Assume that $2K_X$ is Cartier.
Then $(Z\ni o)$ is the quotient of the isolated hypersurface singularity
\[
\{x_1x_3+x_2 \phi(x_2^2,x_4)=0\}\ni 0
\]
by the $\mumu_2$-action given by the weights $(1,1,0,0)$. The contraction $f$ \textup(resp. $f^+$\textup) is the quotient of the blowup of the plane $\{x_2=x_3=0\}$ \textup(resp. $\{x_1=x_2=0\}$\textup) by $\mumu_2$.
In particular, $X$ contains a unique point of index $2$ and the central fiber $C$ is irreducible. The variety $X^+$ is Gorenstein.
\end{theorem}

Similar description is known for divisorial extremal curve germs of index $2$  \cite[\S~4]{KM:92}.

\subsubsection*{First properties}
Let $(X\supset C)$ be an extremal curve germ and let $f: 
(X\supset C)\to (Z\ni o)$ be the corresponding contraction. 
For any connected subcurve $C'\subset C$ the germ $(X\supset C')$ is again an 
extremal curve germ. If moreover $C'\subsetneqq C$, then $(X\supset C')$ is 
birational. 
By the Kawamata-Viehweg vanishing theorem 
\begin{equation}
\label{eq:R1O}
R^1f_*\OOO_X=0
\end{equation} 
(see e.g. \cite{KMM}). As a consequence one has
$\p(C')\le 0$ for any subcurve $C'\subset C$.
In particular, $C=\cup C_i$ is a ``tree'' of smooth rational curves.
Furthermore, 
\begin{equation}
\label{eq:Pic}
\Pic(X)\simeq H^2(X,\ZZ)\simeq \ZZ^{\oplus n},
\end{equation} 
where $n$ the number of irreducible components of $C$.
For more delicate properties of extremal curve germs one needs to know the cohomology 
of the dualizing sheaf, see~\cite{Mori:flip, MP:cb1}:

\begin{equation}
\label{eq:R1w}
R^1f_*\upomega_X=
\begin{cases}
0 & \text{if $f$ is birational},
\\
\upomega_Z& \text{if $f$ is $\QQ$-conic bundle and $Z$ is smooth}.
\end{cases}
\end{equation} 

 \begin{definition*}
 An irreducible extremal curve germ $(X\supset C)$ is (locally) 
\emph{imprimitive} at a point $P$ if 
the inverse image of $C$ under the index-one cover $(X'\ni P')\to (X\ni P)$ 
splits.
\end{definition*}

\begin{theorem}[{\cite{Mori:flip,MP:cb1}}]
\label{thm:sing-points}
Let $(X\supset C)$ be an extremal curve germ and let $C_1,\dots, C_n$ be irreducible components of $C$. 
\begin{itemize}
\item 
Each $C_i$ contains at most $3$ singular points of $X$.
\item 
Each $C_i$ contains at most $2$ non-Gorenstein of $X$ and at most $1$ point which is imprimitive for $(X\supset C_i)$.
\item 
If $X$ is Gorenstein at the intersection point $P=C_i\cap C_j$, $C_i\neq C_j$, then $X$ is smooth outside $P$ 
and $(X\supset C)$ is a $\QQ$-conic bundle germ over a smooth base.
\end{itemize}
\end{theorem}

To prove the first assertion, one needs to analyze the conormal sheaf $I_C/I_C^2$ and use the vanishing 
$H^1(\OOO_X/ J)=0$
for any $J\subset \OOO_X$ with $\Supp(\OOO_X/J)=C$
(see \cite{Mori:flip, MP:1pt}). 
For the second assertion one can use topological arguments based on \eqref{eq:pi1} (see \cite{MP:1pt}).
For the last assertion we refer to \cite[1.15]{Mori:flip}, \cite[4.2]{Kollar:real3}, and~\cite[4.7.6]{MP:1pt}

The techniques applied in the proof of the above proposition allow to obtain much stronger results.
In particular, it allows to classify all the possibilities for the local 
behavior of an irreducible germ $(X\supset C)$ near a singular point $P$ \cite{Mori:flip}.
Thus according to \cite{Mori:flip} and \cite{MP:cb1} 
the triple $(X\supset C\ni P)$ belongs to one of the following types:
\[
\mtypec{IA},\ \mtypec{IC},\ \mtypec{IIA},\ \mtypec{IIB},\
\mtypec{IA^\vee},\ \mtypec{II^\vee},\ \mtypec{ID^\vee},\ \mtypec{IE^\vee},\ \mtypec{III}.
\]
Here the symbol ${}^\vee$ means that $(X\supset C\ni P)$ is locally imprimitive, the symbol \type{II} means that $(X\ni P)$ is a terminal point of exceptional type \type{cAx/4} (see Proposition~\ref{prop:local-ge}), and \type{III} means that $(X\ni P)$ is an (isolated) \type{cDV}-point.

For example, a triple $(X\supset C\ni P)$ is of type \typec{IC} if 
there are analytic isomorphisms 
\begin{equation*} 
(X\ni P)\simeq\CC^3_{y_1,y_2,y_4}/\mumu_m(2,m-2,1), \quad 
C\simeq\{y_1^{m-2}-y_2^2=y_4=0\}/\mumu_m, 
\end{equation*}
where $m$ is odd and $m \ge 5$. 
For definitions other types we
refer the reader to \cite{Mori:flip} and \cite{MP:cb1}.

\subsection{Construction of germs by deformations}
\label{sect-def}
Let $(X\supset C)$ be an extremal curve germ
and let $f: X\to (Z\ni o)$ be the corresponding contraction.
Denote by $|\OOO_Z|$ the infinite dimensional linear system of hyperplane sections passing through $o$ and let $|\OOO_X|:=f^*|\OOO_Z|$. The \emph{general 
hyperplane section} of $(X\supset C)$ is the general member $H\in |\OOO_X|$.
The divisor $H$ contains much more information on the total space than general elephant $D\in |{-}K_X|$. However, the singularities of $H$ typically are more complicated, in particular, $H$ can be non-normal.

The variety $X$ (resp. $Z$) can be viewed as the total space of a one-parameter deformation of 
$H$ (resp. $H_Z:=f(H)$). We are going to reverse this consideration.

\begin{construction*}[see {\cite[\S~11]{KM:92}}, {\cite[\S~1b]{Mori:flip}}]
Suppose we are given a normal surface germ $(H\supset C)$ along a proper curve $C$
and a contraction $f_H : H \to H_Z$ such that $C$ is a fiber and ${-}K_H$ is $f_H$-ample. 
Let $P_1,\dots, P_r \in H$ be all the singular points. Assume also that near each $P_i$ 
there exists
a small one-parameter deformation $\mathfrak{H}_i$ of a neighborhood $H_i$ 
of $P_i$ in $H$ such that the total space $\mathfrak{H}_i$ has a terminal 
singularity at $P_i$. 
The obstruction to globalize deformations
\[
\Psi: \operatorname{Def} (H) \longrightarrow \prod_{P_i\in \Sing(H)} 
\operatorname{Def} (H, P_i)
\]
lies in $R^2 f_* \mathcal{T}_H$, where $\mathcal{T}_H=
\Hom(\Omega_H, \OOO_H )$ is the tangent sheaf of $H$. Since $R^2 f_* \mathcal{T}_H=0$ by the dimension reason, the morphism $\Psi$ is smooth and so there exists a global one-parameter deformation $\mathfrak{H}$ 
of $H$ inducing
a local deformation of $\mathfrak{H}_i$ near $P_i$. 

Then we have a threefold $X:=\mathfrak{H}\supset C$ with 
$H \in |\OOO_X |$ such that locally near $P_i$ it has the desired
structure and one can extend $f_H$ to a contraction $f: X\to Z$
which is birational (resp. a $\QQ$-conic bundle) if $H_Z$ is a surface (resp. a 
curve).
\end{construction*}

\begin{example*}
Consider a rational curve fibration 
$f_{\tilde H} : \tilde H \to H_Z$ over a smooth curve germ $H_Z\ni o$, where $\tilde H$ 
is a smooth surface, such that the fiber over $o$ 
has 
the following weighted dual graph:
\begin{equation}
\label{eq:def-diag}
\vcenter{
\xymatrix@R=3pt{
\overset{-2}{\scriptstyle{\square}}\ar@{-}[r]&\overset{-1}{\bullet}\ar@{-}[r]&\overset{-3}{
\circ}\ar@{-}[r]&\overset{-2}\circ\ar@{-}[r]&\overset{-3}{\circ}\ar@{-}[r]
&\overset{-1}{\bullet}
\\
&&&\underset{-3}{\circ}\ar@{-}[r]\ar@{-}[u]&\underset{-1}{\bullet}
} }
\end{equation}
Contracting the curves corresponding to the white vertices $\scriptstyle{\square}$ and $\circ$ we obtain a singular surface $H$ and a $K_H$-negative contraction $f_{H} : H \to H_Z$ whose
fiber over $o$ is a curve $C\subset H$ having three irreducible components that correspond to the black vertices $\bullet$. 
The singular locus of $H$ consists of a Du Val point $P_0\in H$ of type \type{A_1} and 
a log canonical singularity $P\in H$ 
whose dual graph is formed by the white circle vertices $\circ$.
Both $P_0$ and $P$ have a 1-parameter $\QQ$-Gorenstein smoothings \cite[Computation~6.7.1]{KM:92}.
Applying the above construction to $H\supset C$ we obtain 
an example of a $\QQ$-conic bundle 
contraction $f:(X\supset C)\to (Z\ni o)$ with a unique non-Gorenstein point which is of type \type{cD/3}. 
If we remove the $(-2)$-curve corresponding to $\scriptstyle{\square}$ on the left hand side of the graph \eqref{eq:def-diag}, we get a birational contraction of surfaces 
$f_{H}' : H' \to H_Z'$. Applying the same construction to $H'\supset C$ we obtain 
an example of a divisorial contraction.
Similarly, removing further one of the $(-1)$-curves we get a flip.
\end{example*}

\section{Extremal curve germs: general elephant}

\begin{theorem}[Mori {\cite{Mori:flip}}, Koll\'ar-Mori {\cite{KM:92}}, Mori-Prokhorov {\cite{MP:cb3}}]
\label{thm:ge}
Let $(X\supset C)$ be an irreducible extremal curve germ. Then the 
general member $D\in |{-}K_X|$ has only Du Val singularities. 
\end{theorem}
 
The existence of a Du Val elephant for extremal curve germs with reducible 
central fiber is not known at the moment. See Theorem~\ref{thm:ge-red} below
for partial results in this direction.

\par\smallskip\noindent
\emph{Comment on the proof.}
Essentially, there are three methods to find a good elephant $D\in |{-}K_X|$. We outline them below.

\subsubsection{}
As in Proposition~\ref{prop:local-ge}, 
near each non-Gorenstein point $P_i\in X$ take a local general elephant
$D_i\in |{-}K_{(X\ni P_i)}|$.
Since $D_i$ is general, we have $D_i\cap C=\{P_i\}$. 
Then we can regard $D:=\sum D_i$ as a Weil divisor on $X$.
By the construction, $K_X+D$ is a Cartier divisor near each $P_i$, hence it is Cartier everywhere. In some cases it is possible to compute the intersection numbers $D_i\cdot C$ and show that $D\cdot C<1$.
Then we may assume that $K_X+D\sim 0$ by~\eqref{eq:Pic} and so $D$ is a Du Val anticanonical divisor. For example, 
this method works for extremal curve germs described in Theorems~\ref{thm:cb-index2} and~\ref{thm:indexdf}, as well as in Example~\ref{ex:toroidal} below. 

\subsubsection{}
\label{subsect-2}
In some cases, the above approach does not work but it allows to show the existence of 
a klt $2$-complement $S\in |{-}2K_X|$ such that $\dim (D\cap C)=0$.
Then one can try to extend a good element from the surface $S$.
The crucial 
fact here is that 
the natural map 
\begin{equation*}
\tau: H^0(X,\OOO_X({-}K_X))\longrightarrow H^0(S,\OOO_S({-}K_X))=\upomega_{S}
\end{equation*}
is surjective if $(X, C)$ is birational 
and surjective modulo $\Omega_{S}^2$ if $(X, C)$ is a $\QQ$-conic bundle.
This immediately follows from \eqref{eq:R1w}.
Details can found in \cite[\S~2]{KM:92} and \cite{MP:cb3}.

\subsubsection{}
Finally, in the most complicated cases none of the above methods work.
Then one needs more subtle techniques which requires detailed analysis of singularities and infinitesimal structure of $X$ along $C$ \cite[\S\S~8-9]{Mori:flip}. Then, roughly speaking, the good section $D\in |{-}K_X|$ is recovered as the formal Weil divisor $\underset{\longleftarrow}\lim\, C_n$ of the completion $X^{{\wedge}}$ of $X$ along $C$, where $C_n$ are subschemes with support $C$ constructed by using certain inductive procedure \cite[\S~9]{Mori:flip}.

As a consequence of Theorem~\ref{thm:ge}, in the $\QQ$-conic bundle case 
one obtains the following fact which proves Iskovskikh's conjecture \cite{Isk96:conic-e}. 

\begin{corollary}
\label{cor:base}
Let $(X\supset C)$ be a $\QQ$-conic bundle germ over $(Z\ni o)$, where $C$ can be reducible. 
Then $(Z\ni o)$ is a Du Val singularity of type \type{A_n} (or smooth).
\end{corollary}

This result is very useful for applications to the rationality problem of 
three-dimensional varieties having conic bundle structure \cite{Isk96:conic-e, P:rat-cb:e} and some problems of biregular geometry \cite{P:degQ-Fano-e, P:2010:QFano}.

It turns out that the structure of $\QQ$-conic bundle germs over a singular base $(Z\ni o)$ is much simpler and short than other ones. In fact these germs can be exhibited as certain quotients of
$\QQ$-conic bundles of index $\le 2$ (see Theorem~\ref{thm:cb-index2}). A complete classification of such germs was obtained in \cite{MP:cb1, MP:cb2}. Here is a typical example.

\begin{example}
\label{ex:toroidal}
Let the group $\mumu_{n}$ act on on $\CC^2
_{u,v}$ and $\PP^1 _{x,y} \times\CC^2 _{u,v}$ via
\begin{equation*}
(x {:}\, y;\, u,\, v)\longmapsto \big(x {:}\, \zeta^ay;\, \zeta u,\,
\zeta^{-1}v\big),
\end{equation*}
where $\zeta=\zeta_n =\exp(2\pi i/n)$ and $\gcd(n,a)=1$.
Then the projection 
\[
 f: X=\big(\PP^1 \times\CC^2\big)/\mumu_{n}\longrightarrow Z=\CC^2
/\mumu_{n}
\]
is a $\QQ$-conic bundle. 
The variety $X$ has 
exactly two singular points which are terminal cyclic
quotients of type $\frac{1}{n}(1,-1,\pm a)$. The surface $Z$ has
at $0$ a Du Val point of type \type{A_{n-1}}. 
\end{example}
McKernan proposed a natural higher-dimensional analogue of Corollary~\ref{cor:base}:

\begin{conjecture}
Let $f: X\to Z$ be a $K$-negative contraction such that $\uprho(X/Z)=1$ and $X$ is 
$\varepsilon$-lc, that is, all the coefficients in \eqref{eq:discr} satisfy $a_i\ge -1+\varepsilon$. 
Then $Z$ is $\delta$-lc, where $\delta$ depends 
on $\varepsilon$ and the dimension.
\end{conjecture}
A deeper version of this conjecture which generalizes  Theorem~\ref{thm:DP-mult} and uses the notion was proposed by Shokurov. 
He also suggested that the optimal value of $\delta$, in the case where singularities of $X$ are canonical and $f$ has one-dimensional fibers, equals $1/2$. 
Recently, this was proved by 
J.~ Han, C.~Jiang, and Y.~Luo \cite{HanJiangnLuo:shokurovs}.

Once we have a Du Val general elephants,
all extremal curve germs
can be divided into two large classes
which will be discussed
separately in the next two sections.

\begin{definition}
Let $(X\supset C)$ be an extremal curve germ
and let $f: X\to (Z\ni o)$ be the corresponding contraction. Assume that the general member 
$D\in |{-}K_X|$ is Du Val. Consider the 
Stein factorization:
\[
f_D: D \longrightarrow D' \longrightarrow f(D)\qquad \qquad(\text{put $D'=f(D)$ if 
$f$ is birational}).
\]
Then the germ $(X\supset C)$ is said to be \emph{semistable} if 
$D'$ has only (Du Val) singularities of type \type{A_n}.
Otherwise, $(X\supset C)$ is called \emph{exceptional}. 
\end{definition}

\section{Semistable germs}
Let $(X\supset C)$ be an irreducible extremal curve germ.
By Theorem~\ref{thm:ge} the general member $D\in |{-}K_X|$ is Du Val.
In this section we assume that $(X\supset C)$ is semistable. 
Excluding simple cases, we assume also that $X$ is not Gorenstein \cite{Cutkosky:contr} and $(X\supset C)$ is not a $\QQ$-conic bundle germ over a singular base \cite{MP:cb1,MP:cb2}.
According to Theorem~\ref{thm:sing-points} the threefold $X$ has at most two non-Gorenstein points. 
Thus the following case division is natural:

\begin{itemize}
\item 
[Case \typec{k1A}:]
the set of non-Gorenstein points consists of a single point $P$;

\item 
[Case \typec{k2A}:]
the set of non-Gorenstein points consists of exactly two points $P_1$, $P_2$.
\end{itemize}

\begin{proposition}
In the above hypothesis, for the general member $H\in |\OOO_X|$ 
the pair $(X, H+D)$ is lc. If moreover $D\supset C$, 
then 
$H$ is normal and has only cyclic quotient singularities.
In this case the
singularities of $H$ are of type~\type{T}. 
\end{proposition}

The proof uses the 
inversion of adjunction \cite{Shokurov:flips+err} to extend 
a general hyperplane section from $D$ to $X$ (see \cite[Proposition~2.6]{MP:IA}).

For an extremal curve germ of type \typec{k2A} any member $D\in |{-}K_X|$ contains $C$ \cite{KM:92}.
Hence the general hyperplane section $H\in |\OOO_X|$ has only \type{T}-singularities
and $X$ can be restored as a one-parameter deformation space of $H$. 
In this case
 $X$ has no singularities other than $P_1,\, P_2$.
Moreover, $(X\supset C)$ cannot be a $\QQ$-conic bundle germ \cite{MP:cb1,MP:cb3}. The birational germs of type \typec{k2A} were completely described by Mori \cite{Mori:ss}. He gave an explicit algorithm for computing divisorial contractions and flips in this case.

The structure of extremal curve germs of type \typec{k1A} is more complicated.
They were studied in \cite{MP:IA}. 
In particular, the general hyperplane section $H\in |\OOO_X|$ was computed.
However \cite{MP:IA} does not provide a good description of the infinitesimal structure
of $X$ along $C$ neither an algorithm similar to \cite{Mori:ss}.
This was done only in a special situation in \cite{HTU}.
Note that in the case \typec{k1A} a general member $H\in |\OOO_X|$ can be non-normal.

\begin{examples*}
Similar to the example in \S~\ref{sect-def}, consider a surface germ $H\supset C\simeq \PP^1$ whose dual graph has the following graph of the minimal 
resolution:
\[
\xymatrix{\overset {-1} \bullet\ar@{-}[r]&\overset{-7}\circ\ar@{-}[r]&
\overset{-2}\circ\ar@{-}[r]&\overset{-2}\circ\ar@{-}[r]&\overset{-2}\circ} 
\]
where $\bullet$ is a $(-1)$-curve.
The chain formed by white circle vertices $\circ$ corresponds to a \type{T}-singularity of type 
$\frac1{25}(1,4)$.
The whole configuration can be contracted to a cyclic quotient singularity $H_Z\ni o$ of type
 $\frac1{21}(1,16)$. Since this is not a 
\type{T}-singularity, the induced threefold contraction must be flipping.
\end{examples*}

\section{Exceptional curve germs}

In this section we assume that $(X\supset C)$ is an exceptional irreducible extremal curve germ. 
As in the previous section we also assume that $X$ is not Gorenstein and $(X\supset C)$ is not a $\QQ$-conic bundle germ over a singular base.
According to the classification \cite{Mori:flip,KM:92,MP:cb3}
the germ $(X\supset C)$ belongs to one of
following types:

\begin{itemize}
\item
$X$ has a unique non-Gorenstein point $P$ which is of type 
\type{cD/2}, \type{cAx/2}, \type{cE/2}, or 
\type{cD/3} and $(X\supset C)$ is of type \typec{IA} at $P$;

\item
$X$ has a unique non-Gorenstein point $P$ which is of 
exceptional type \type{cAx/4}
and $(X\supset C)$ is of type 
\typec{IIA}, \typec{II^\vee}, or \typec{IIB} at $P$;

\item 
$X$ has a unique singular point $P$ which is a cyclic quotient singularity 
of index $m\ge 5$ \textup(odd\textup) and $(X\supset C)$ is of type \typec{IC} at $P$;

\item 
$X$ has two singular points of indices $m\ge 3$ \textup(odd\textup) and $2$, then $(X\supset C)$ is said to be of type \typec{kAD};
\item 
$X$ has three singular points of indices $m\ge 3$ \textup(odd\textup), $2$ and $1$, then $(X\supset C)$ is said to be of type \typec{k3A}.
\end{itemize}
In each case 
the general elephant is completely described in terms of its minimal resolution:

\begin{theorem}[{\cite{KM:92,MP:cb3}}]
\label{theorem:ge}
In the above hypothesis assume that the general element
$D\in |{-}K_X|$ contains $C$. 
Then the dual graph of
$(D\supset C)$ is one of the following, where white vertices $\circ$ denote $(-2)$-curves on the minimal resolution of $D$ and the black vertex $\bullet$ corresponds to the proper transform of $C$:
\par\smallskip\noindent
$
\xymatrix@R=-1pt@C=17pt{
\mathrm{{(IC)}}&{\underbrace{\circ -\cdots - \circ}_{m-3\ge 2}} 
\ar@{-}[r]&\circ\ar@{-}[r]\ar@{-}[d]&\circ
\\
&&\bullet&
}\qquad\qquad
\xymatrix@R=4pt@C=13pt{
\mathrm{{(IIB)}}&&&\circ\ar@{-}[d]
\\
&\circ\ar@{-}[r]&\circ\ar@{-}[r]&\circ\ar@{-}[r]&\circ\ar@{-}[r]&\bullet 
} 
$

\par\smallskip\noindent
$\mathrm{{(kAD)}}$ $
\xymatrix@R=4pt@C=11pt{
&&&&&\circ\ar@{-}[d]
\\
&{{\circ -\cdots - \circ}} \ar@{-}[r]& {\bullet}\ar@{-}[r]& \circ\ar@{-}[r] 
&\cdots\ar@{-}[r] &\circ&\circ\ar@{-}[l]
} 
$
\par\smallskip\noindent
$\mathrm{(k3A)}$\hspace{20pt}
$
\xymatrix@R=4pt@C=17pt{
&\circ
\\
{{\circ -\cdots - \circ}} \ar@{-}[r]&\bullet \ar@{-}[u]\ar@{-}[r]&\circ
} 
$
\end{theorem}

Exceptional irreducible 
extremal curve germs are studied are well (see \cite{KM:92}, \cite{MP:1pt}, and references therein).
For flipping ones the general hyperplane section $H\in |\OOO_X|$ is normal and has only rational singularities.
It is computed in \cite{KM:92} and the flip is reconstructed as a
one-parameter deformation space of $H$.
For divisorial and $\QQ$-conic bundle germs, the situation is more complicated. 
Then the general hyperplane section $H$ can be non-normal (see e.g. \cite{MP:IIAb}). Nevertheless in almost all cases, except for types \typec{kAD}and \typec{k3A}, there is a description of $H\in |\OOO_X|$ and infinitesimal structure of these germs.
For convenience of references in the table below we collect 
the known information on the exceptional irreducible extremal curve germs.

\par\smallskip\noindent
\begin{tabular}{m{0.2\textwidth}ll} 
type & $(X,C)$ & references
\\[-0.7em]
\multicolumn{3}{c}{\rule{0.97\textwidth}{.3mm}}
\\[-0.4em]
index $2$ germs &
divisorial, $\QQ$-conic bundle & \cite[\S~4]{KM:92}, \cite[\S~12]{MP:cb1}, \cite[\S~7]{MP:IA}
\\
\type{cD/3}
&
flip, divisorial &\cite[\S~6]{KM:92}, \cite[\S~4]{MP:IA}
\\
\typec{IC} &
flip, $\QQ$-conic bundle (only for $m=5$) &
\cite[\S~8]{KM:92}, \cite{MP:IC-IIB}
\\
\typec{IIA} &
flip, divisorial, $\QQ$-conic bundle &
\cite[\S~7]{KM:92}, \cite{MP:IIAb,MP:IIAa} 
\\
\typec{IIB} &
divisorial, $\QQ$-conic bundle &
\cite{MP:IC-IIB} 
\\
\typec{II^\vee} &
divisorial, $\QQ$-conic bundle &
\cite[4.11.2]{KM:92}, \cite{MP:cb1} 
\\
\typec{kAD} & flip, divisorial, 
$\QQ$-conic bundle &
\cite[\S~9]{KM:92}, \cite{Mori:errKM, MP:cb1,MP:cb3}
\\
\typec{k3A} &
divisorial, $\QQ$-conic bundle &
\cite[\S~5]{KM:92}, 
\cite{MP:cb1, MP:cb3} 
\end{tabular}
\par\smallskip\noindent
Detailed analysis of the local structure of exceptional extremal curve germs
allows to extend the result of 
Theorem~\ref{thm:ge} to the case of reducible central fiber
containing an exceptional component:

\begin{theorem}[Mori-Prokhorov {\cite{MP:ge}}]
\label{thm:ge-red}
Let $(X\supset C)$ be an extremal curve germ such that
$C$ is reducible and satisfies the following condition:
\begin{enumerate}
\item[(*)] 
each component $C_i\subset C$ contains at most one point of index $>2$.
\end{enumerate}
Then the general member $D\in |{-}K_X|$ has only Du Val singularities.
Moreover, for each irreducible component $C_i\subset C$ with two non-Gorenstein 
points or of types \typec{IC} or \typec{IIB}, the dual graph of $(D,C_i)$
has the same form as the irreducible extremal curve germ $(X\supset C_i)$.
\end{theorem}

The proof uses the extension techniques of sections of $|{-}K_X|$ from a good member $S\in |{-}2K_X|$ (see~\ref{subsect-2}).

\section{$\QQ$-Fano threefolds}
In arbitrary dimension $\QQ$-Fano threefolds are bounded, i.e. they are contained in fibers of a morphism of schemes of finite type.
This is a consequence of the much more general fact \cite{Birkar-BAB}.
In dimension $3$ there are effective results based on the orbifold Riemann-Roch formula \eqref{eq:RR} and Bogomolov-Miyaoka inequality applied to the restriction of the reflexive sheaf $(\Omega^1_X)^{\vee\vee}$ to a sufficiently general hyperplane section \cite{Kawamata:bF}. In particular, combining \eqref{eq:RR} with Serre duality we obtain
\[
\chi (\OOO_X) = \frac1{24} \left({-}K_X\cdot \mathrm{c}_2(X) +\sum_P \left(m_P-\frac1{m_p}\right) \right)
\]
where $m_P$ is the index of a virtual quotient singularity of $X$ \cite{Reid:YPG}. Since $X$ is $\QQ$-Fano, by Kawamata-Viehweg vanishing theorem \cite{KMM} one has 
$\chi (\OOO_X)=1$. Arguments based on Bogomolov-Miyaoka inequality shows that 
${-}K_X\cdot \mathrm{c}_2(X)$ is positive (see \cite{Kawamata:bF}). This gives an effective bound of 
the indices of singularities of $X$. Similarly one can get an upper bound of the anticanonical degree ${-}K_X^3$. Moreover, analyzing the methods of \cite{Kawamata:bF} it is possible to enumerate Hilbert series of all $\QQ$-Fano threefolds. This information is collected in \cite{GRD} in a form of a huge computer database of possible ``candidates''. It was extensively explored by many authors, basically, to obtain lists of examples representing $\QQ$-Fano threefolds as subvarieties of small codimension in a weighted projective space
(see e.g. \cite{Fletcher:wci, Brown-Kerber-Reid:codim4}) and references therein).

\begin{examples*}
\begin{itemize}
 \item 
There are 130 (resp. 125) families of $\QQ$-Fano threefolds that are representable as hypersurfaces (resp. codimension two complete intersections) in weighted projective spaces \cite{Fletcher:wci,GRD}. 
 \item 
Toric $\QQ$-Fano threefolds are exactly weighted projective spaces 
$\PP(3, 4, 5, 7)$, 
$\PP(2, 3, 5, 7)$, 
$\PP(1, 3, 4, 5)$, 
$\PP(1, 2, 3, 5)$, 
$\PP(1,1, 2, 3)$, 
$\PP(1,1,1, 2)$, 
$\PP^3=\PP(1,1,1, 1)$, 
and the quotient of $\PP^3$ by $\mumu_5$ that acts diagonally with weights $(1,2,3,4)$ \cite{GRD}.
\end{itemize}
\end{examples*}

Although the classification is very far from completion, there are several
systematic results. For example, the optimal upper bound of the degree ${-}K_X^3$ of $\QQ$-Fano threefolds was obtained in \cite{P:degQ-Fano-e}. If $X$ is singular, it is equal to $125/2$ and achieved for the weighted projective space $\PP(1,1,1,2)$.
The lower bound of the degree equals $1/330$ \cite{Chen-Chen:Fano} and is achieved for a hypersurface of degree $66$ in $\PP(1,5,6,22,33)$.
It is known that, under certain conditions, General Elephant Conjecture~\ref{conj:ge} holds 
for $\QQ$-Fano threefolds modulo deformations \cite{Sano:eleph}.

\begin{ack}
The author would like to thank Professors Shigefumi Mori and Vyacheslav  Shokurov  for helpful comments on the original version of this paper.
\end{ack}

\begin{funding}
This work was performed at the Steklov International Mathematical Center and supported by the Ministry of Science and Higher Education of the Russian Federation (agreement no. 075-15-2019-1614).
\end{funding}

\def\cprime{$'$}

\end{document}